
\magnification=\magstep1
 \hsize=12.5cm
 \vsize=18cm
 \hoffset=1cm
 \voffset=2cm

\def\hf{{\textstyle{1\over2}}}
\def\txt#1{{\textstyle{#1}}}
\def\no{\noindent}
\def\a{\alpha} \def\b{\beta}
\def\d{{\,\rm d}}\def\D{\delta}
\def\e{\varepsilon}

\def\k{\kappa}
\def\s{\sigma}
\def\t{\theta}
\def\={\;=\;}

\def\zt{\zeta(\hf+it)}

\def\DJ{\leavevmode\setbox0=\hbox{D}\kern0pt\rlap
 {\kern.04em\raise.188\ht0\hbox{-}}D}

\font\tenmsb=msbm10
\font\sevenmsb=msbm7
\font\fivemsb=msbm5
\newfam\msbfam
\textfont\msbfam=\tenmsb
\scriptfont\msbfam=\sevenmsb
\scriptscriptfont\msbfam=\fivemsb
\def\Bbb#1{{\fam\msbfam #1}}

\def \NN {\Bbb N}

\def \ZZ {\Bbb Z}

\font\aa=cmss12
\font\bb=cmcsc10

\font\hh=cmr10 at 9pt

\centerline{\bb  ON SOME RECENT RESULTS IN THE THEORY OF THE
ZETA-FUNCTION}
\medskip
\centerline{\aa Aleksandar Ivi\'c}
\medskip
\centerline{\sevenbf in `Contemporary Mathematics - 125 years of Faculty
of Mathematics' ed. N. Bokan}
\centerline{\sevenbf University of Belgrade, Belgrade, 2000, pp. 83-92}
\medskip
\centerline{\bb 1. Introduction}
\medskip\no
This review article is devoted to
the Riemann zeta-function $\zeta(s)$,  defined for $\Re\,s > 1$ as
$$
\zeta(s) = \sum_{n=1}^\infty n^{-s} =
\prod_{p\,:\,{\rm{prime}}} (1 - p^{-s})^{-1}, \leqno{(1.1)}
$$
and otherwise by analytic continuation. It admits meromorphic
continuation to the whole complex plane, its only singularity
being the simple pole $s = 1$ with residue 1. For general
information on $\zeta(s)$ the reader is referred to the
monographs [13], [16], [21], [28] and [51]. From the functional equation
$$
\zeta(s) = \chi(s)\zeta(1 - s), \quad\chi(s) = 2^s\pi^{s-1}
\sin\bigl({{\pi s}\over 2}\bigr)\Gamma(1 - s), \leqno{(1.2)}
$$
which is valid for any complex $s$, it follows that $\zeta(s)$
has zeros at $s = -2,-4,\ldots$~. These zeros are 
called the ``trivial'' zeros of $\zeta(s)$, to distinguish them
from the complex zeros of $\zeta(s)$. The zeta-function has
also an infinity of complex zeros.
It is well-known that all complex zeros of  $\zeta(s)$  lie
in the so-called ``critical strip'' $0 < \sigma = \Re\, s < 1$, and
if $N(T)$ denotes the number of zeros $\rho = \beta + i\gamma$
($\beta, \gamma$ real) of  $\zeta(s)$ for which $0 < \gamma \le T$,
then ($f = O(g)$ and $f \ll g$ both mean that $|f(x)| \le Cg(x)$
for some $C > 0$ and $x \ge x_0$)
$$
N(T) = {T\over {2\pi}}\log\bigl({T\over {2\pi}}\bigr) - {T\over {2\pi}}
+ {7\over 8} + S(T) + O({1\over T}) \leqno{(1.3)}
$$
with
$$
S(T) = {1\over \pi}\arg \zeta(\txt{1\over 2} + iT) = O(\log T). \leqno{(1.4)}
$$
This is the so-called Riemann--von Mangoldt formula. {\it The 
Riemann hypothesis} (henceforth RH for short) is the
conjecture, stated by B. Riemann in his epoch-making memoir [49], that
{\it very likely -- sehr wahrscheinlich -- all complex zeros of  $\zeta(s)$
have real parts equal to $1\over2$}. For this reason the line
$\sigma = \txt{1\over2}$ is called the ``critical line'' in the theory of 
$\zeta(s)$. The RH is undoubtedly one of the most celebrated and difficult
open problems in whole Mathematics. Its proof  would
have very important consequences in multiplicative number theory,
especially in problems involving the distribution of primes, since
the defining relation (1.1) shows the intrinsic connection between
$\zeta(s)$ and primes. It would also
very likely lead to generalizations to many other zeta-functions
(Dirichlet series) sharing similar properties with  $\zeta(s)$.
Despite much evidence in favour of the RH, there are also some reasons
to be skeptical about its truth -- see, for example, [24], [27] and [30].
The aim of this paper is to present briefly some recent results in the
theory of $\zeta(s)$. The choice of subjects is motivated by the 
limited length of this text and by the author's personal research
interests. An extensive bibliography will aid the interested reader.

\medskip
\centerline{\bb 2. Zeros on the critical line}
\medskip
\no The distribution of the zeros $\rho = \beta + i\gamma$ of $\zeta(s)$
on the critical line $\s = \hf$
is a fundamental problem in the theory of $\zeta(s)$. Despite enormous
efforts the RH has been so far neither proved nor disproved. The
smallest zeros of $\zeta(s)$
(in absolute value) are ${1 \over 2} \pm 14.134725\ldots i$.
Large scale computations of zeros of $\zeta(s)$ have been
carried out in recent times by  the aid of computers (see [44]--[46], [48]).
 Suffice to say that it is known today that the first 1.5 billion
complex zeros of  $\zeta(s)$ in the upper half-plane are simple
and do have real parts equal to $\txt{1\over2}$, as predicted by the RH.
Moreover, many large blocks of zeros of much greater height
have been thoroughly investigated, and all known zeros satisfy
the RH. The RH is, however, not the only topic connected with the
zeros on  the critical line.  Let $N_0(T)$ denote the number of
zeros  of $\zeta(s)$ of the form $\hf + it, 0 < t \le T$. The RH is
in fact the statement that $N_0(T) = N(T)$ for $T > 0$. A fundamental
result of A. Selberg [50] from 1942 states that 
$$N_0(T) > CN(T)\leqno(2.1)$$ 
for $T \ge T_0$
and some positive constant $C > 0$. In other words, a positive proportion
of all complex zeros of $\zeta(s)$ lies on the critical line.
A substantial advancement in this topic was made by N. Levinson [39]  in
1974, who proved that one can take $C = 1/3$ in (2.1). By refining 
Levinson's techniques further improvements have been made, and today it
is known  that one can take $C = 2/5$  in (2.1), but there is no
chance that by existing methods one can reach the value $C = 1$.

\medskip
Another interesting problem is the problem of the gap between 
consecutive zeros on the critical line.
Let us denote by $0 < \gamma_1 \le \gamma_2 \le \ldots$ the
positive zeros of $\zeta(\hf + it)$ with multiplicities counted.
 If the RH is true, then it is known
(see [51]) that the function $S(T)$, defined by (1.3) and (1.4),
satisfies
$$
S(T) = O\Bigl({\log T \over \log \log T} \Bigr). \leqno{(2.2)}
$$
This seemingly small improvement over (1.4) is significant:
If (2.2) holds, then
from (1.3) one infers  that $N(T + H) - N(T) > 0$ for
$H = C/\log \log T$ with a suitable $C > 0$ and $T \ge T_0$.
Consequently we have, assuming the RH, the bound
$$
\gamma_{n+1} - \gamma_n \, \ll \, {1 \over \log \log \gamma_n} 
\leqno{(2.3)}
$$
for the gap between consecutive zeros
on the critical line. The bound (2.3) is certainly out of reach
at present. For some unconditional results
on $\gamma_{n+1} - \gamma_n$, see [17]--[19] and [33]. For example,
it was proved by the author [16], [17] that one has unconditionally
$$
\gamma_{n+1} - \gamma_n \, \ll \, \gamma_n^{\t+\e},
\qquad \t \= {\kappa+\lambda\over4(\kappa+\lambda)+2},
$$
where $\,(\kappa,\,\lambda)\,$ is a so-called ``exponent pair"
from the theory of exponential sums (see [8], [16] for definition),
and $\e$ denotes arbitrarily small positive constants. 
With known results on exponent pairs one obtains then as the strongest known
result $\t \le 0.155945\ldots\;$. Another unconditional result
of the author (see [19]) states that
$$
\sum_{\gamma_n\le T}(\gamma_{n+1} - \gamma_n)^3 \ll T\log^6T.
$$

\medskip
\centerline{\bb 3. Zero-density results}
\medskip\no
In the absence of a proof of RH, in practice one has to assume
the existence of the zeros of $\zeta(s)$ in $\hf < \s <1$ and show
that the influence of these zeros is not large. One defines the
zero-counting function
$$
N(\s,T) \= \sum_{\zeta(\b+i\gamma)=0,\b\ge\s,|\gamma|\le T}1
$$
and tries to obtain upper bounds for $N(\s,T)$. By the symmetry of zeros
it is  sufficient to assume that $\s \ge \hf$, and one also often
takes into  account the so-called ``zero-free region" of $\zeta(s)$.
This is due to I.M. Vinogradov (see Ch. 6 of [16]) and says that
$$
\b \;\le\;1 - C(\log\gamma)^{-2/3}(\log\log\gamma)^{-1/3}
\qquad(\zeta(\b+i\gamma) =0,\;\gamma \ge \gamma_0 > 0,\; C > 0).\leqno(3.1)
$$
This result was obtained in 1958, and is one of the oldest results in
zeta-function theory which has not been subsequently improved.

\medskip In obtaining upper bounds for $N(\s,T)$ one can use several
techniques (see Ch. 11 of [16] for a comprehensive account).
In the range $\hf < \s \le {3\over4}$ the best bound is due to A.E. Ingham
[10]. This is
$$
N(\s,T) \;\ll\; T^{3(1-\s)/(2-\s)}\log^5T,\leqno(3.2)
$$
while in the range ${3\over4} \le \s \le 1$ one has M.N. Huxley's result [6]
that
$$
N(\s,T) \;\ll\; T^{3(1-\s)/(3\s-1)}\log^{44}T.\leqno(3.3)
$$
When combined, (3.2) and (3.3) yield
$$
N(\s,T) \;\ll\; T^{12(1-\s)/5}\log^{44}T\qquad(\hf \le \s < 1).\leqno(3.4)
$$
From (3.4) one can deduce the bound
$$
p_{n+1} - p_n \ll p_n^{7/12}\log^Cp_n,
$$
where $p_n$ is the $n$--th prime (see Ch. 12 of [16]). By 
combining analytic and sieve techniques
various authors have reduced the exponent ``7/12", and the current record
holders are Baker and Harman [2], who proved that
$$
p_{n+1} - p_n \ll p_n^{0.535}.
$$
It is interesting that the so-called ``density hypothesis" 
$$
N(\s,T) \;\ll\; T^{2-2\s+\e} \qquad(\hf \le \s \le 1)\leqno(3.5)
$$
gives the bound
$$
p_{n+1} - p_n \ll p_n^{{1\over2}+\e},
$$
while the much stronger RH gives only a slightly better result, namely
$$
p_{n+1} - p_n \ll p_n^{{1\over2}}\log p_n,
$$
which is still insufficient to prove the old conjecture: {\it between
every two squares there is always a prime}.

\medskip The bound (3.5) is known to hold for $\s  \ge {11\over14}$, which
was proved by M. Jutila [38] in 1977. For various values of  $\s$ in the
range ${3\over4} < \s < 1$ sharper results than (3.3) or 
(3.4) have been proved. For example, the author [11], [14]--[16] obtained
the bounds
$$\eqalign{&
N(\s,T) \;\ll\; T^{3(1-\s)/(2\s)+\e}\qquad({\txt{3831\over4791}} =
0.799624\ldots \le \s \le 1),
\cr&
N(\s,T) \;\ll\; T^{9(1-\s)/(7\s-1)+\e}\qquad({\txt{41\over53}} =
0.773585\ldots \le \s \le 1),
\cr&
N(\s,T) \;\ll\; T^{6(1-\s)/(5\s-1)+\e}\qquad({\txt{13\over17}} 
= 0.764705\ldots \le \s \le 1).
\cr}$$
\medskip
\centerline{\bb 4. The mean square}
\medskip\no
Power moments, and especially the mean square formulas, play
an important r\^ole in the theory of $\zeta(s)$. In fact, zero-density
estimates discussed in the previous section depend heavily on power
moments of $\zeta(s)$, as discussed extensively in Chapter 11 of [16].
As far as mean square formulas are concerned, one can distinguish
between the cases $\s = \hf$ and $\hf < \s \le 1$. Here we shall
briefly discuss only the former case, and for the latter we refer
the reader to [34], [41]. One has the asymptotic formula
$$
\int_0^T|\zt|^2\d t \= T\log\left({T\over2\pi}\right) + (2\gamma
- 1)T + E(T),\leqno(4.1)
$$
where $\gamma = 0.577\ldots\,$ is Euler's constant, and 
$E(T)$ is to be considered as the  error term in the
asymptotic formula (4.1). F.V. Atkinson [1] established in 1949 an explicit,
albeit complicated formula for $E(T)$, containing two exponential sums of
length $\asymp T$, plus an error term which is $O(\log^2T)$.
Atkinson's formula has been the starting point for many results
on $E(T)$. It is conjectured that $E(T) \ll T^{1/4+\e}$, and
this bound cannot be proved even if the RH is assumed. The 
best known upper bound for $E(T)$, obtained by intricate estimation of
a certain exponential sum, is due to M.N. Huxley [7]. This is
$$
E(T) \;\ll\; T^{72/227}(\log T)^{679/227}, \quad {72\over227} = 
0.3171806\ldots\;.
$$
In the other direction, J.L. Hafner and the author [3], [21] proved 
in 1987 that there exist absolute constants $A, B > 0$ such  that
$$
E(T) = \Omega_+\left\{(T\log T)^{1/4}(\log\log T)^{(3+\log4)/4}
e^{-A\sqrt{\log\log\log T}}\right\}
$$
and
$$
E(T) = \Omega_-\left\{T^{1/4}\exp\left({B(\log\log T)^{1/4}
\over(\log\log\log T)^{3/4}}\right)\right\}.
$$
The omega-symbols are customarily defined as follows:
$f = \Omega(g)$ means that $f = o(g)$ does not hold, 
$f = \Omega_+$ means that $\limsup f/g > 0$, $f = \Omega_-$ means that 
$\liminf f/g < 0$, and 
$f = \Omega_\pm(g)$ means that $\limsup f/g > 0$ and  $\liminf f/g < 0$
both hold. A quantitative $\Omega$--result for $E(T)$ 
was proved by the author [23]: There exist constants $C, D > 0$
such that for $T \ge T_0$ every interval $\,[T,\,T + CT^{1/2}]\,$
contains points $t_1, t_2$ for which
$$
E(t_1) > Dt_1^{1/4},\qquad E(t_2) < -Dt_2^{1/4}.
$$
Numerical calculations pertaining to $E(T)$ have been carried out in [37]
by the author and H. te Riele. Power moments of $E(T)$ were considered
in 1983 by the author  [12], where it was proved that
$$
\int_0^T|E(t)|^A\d t \;\ll\; T^{(A+4)/4+\e}\qquad(0 \le A \le {\txt{35\over
4}}),
$$
which in view of the $\Omega$--results is (up to ``$\e$") optimal,
and supports the conjecture $E(T) \ll T^{1/4+\e}$. Later
in 1992 D.R. Heath-Brown [5] used these bounds to investigate the distribution
function connected with $E(T)$. In the special case $A = 2$ it is known
that
$$
\int_0^T E^2(t)\d t \;=\; {2C\over3}(2\pi)^{-1/2}T^{3/2} + O(T\log^4T),
\qquad C = \sum\limits_{n=1}^\infty d^2(n)n^{-3/2}
\leqno(4.2)
$$
where $d(n)$ is the number of divisors of $n$. The bound for
the error term in (4.2) was obtained independently by E. Preissmann [47]
and the author [21].

In 1992 K.-M. Tsang [52] proved that, for some $\delta > 0$ and explicit
constants $c_1, c_2 > 0$,
$$
\int_0^TE^3(t)\d t \;=\; c_1T^{7\over4} + 
O(T^{{7\over4}-\delta}),
$$ $$
\int_0^TE^4(t)\d t \;=\; c_2T^2 + 
O(T^{2-\delta}).
$$
The author [31] recently improved these bounds to
$O(T^{{47\over28}+\e})$ and $O(T^{{45\over23}+\e})$, respectively.

\medskip
\centerline{\bb 5. The mean fourth power}
\medskip\no
The asymptotic formula for the fourth moment of the Riemann zeta-function
$\zeta(s)$ on the critical line is customarily written as
$$
\int_0^T|\zt|^4\d t \;=\; TP_4(\log T) \;+\;E_2(T),\;
P_4(x) \;=\; \sum_{j=0}^4\,a_jx^j.\leqno(5.1)
$$
A classical result of A.E. Ingham [9] (see Ch. 5 of [16] for a relatively
simple proof) is that $a_4 = 1/(2\pi^2)$ and that the error term $E_2(T)$
in (5.1) satisfies the bound $E_2(T) \ll T\log^3T$. In 1979 D.R. Heath-Brown 
[4] made significant progress in this problem by proving that 
$E_2(T) \ll T^{7/8+\e}$. He also calculated
$$
a_3 \;=\; 2(4\gamma - 1 - \log(2\pi) - 12\zeta'(2)\pi^{-2})\pi^{-2}
$$
and produced more complicated expressions for $a_0, a_1$ and $a_2$ in (5.2). 
For an explicit evaluation of
the $a_j$'s  in (5.1) the reader is referred to the author's work [27]. 
At present the bound
$$
\int_0^T\,|\zt|^k\d t \;\ll_\e\; T^{1+\e}
$$
is not known to hold (see Chs. 7-8 of [16], and [21]) 
for any constant $k > 4$,
which makes the function $E_2(T)$ particularly important in the theory
of mean values of $\zeta(s)$. In recent years, 
due primarily to the application
of powerful methods of spectral theory (see Y. Motohashi's monograph [49] for
a comprehensive account), much advance has been made in connection with
$E_2(T)$.  This
involves primarily exponential sums involving the quantities
$\kappa_j$ and $\a_jH_j^3(\hf)$. Here as usual 
$\,\{\lambda_j = \k_j^2 + {1\over4}\} \,\cup\, \{0\}\,$
is the discrete spectrum of the non-Euclidean Laplacian acting
on $\,SL(2,\ZZ)\,$--automorphic forms, and 
$\a_j = |\rho_j(1)|^2(\cosh\pi\k_j)^{-1}$, where
$\rho_j(1)$ is the first Fourier coefficient of the Maass wave form
corresponding to the eigenvalue $\lambda_j$ to which the Hecke series
$H_j(s)$ is attached. It is conjectured that $E_2(T) \ll T^{1/2+\e}$,
which would imply the (hitherto unproved) bound $\zt \ll t^{1/8+\e}$.
It is known now that
$$
E_2(T) \;=\; O(T^{2/3}\log^{C_1}T),\quad E_2(T) \;=\; \Omega(T^{1/2}),
\leqno(5.2) $$
$$
\int_0^TE_2(t)\d t \;=\; O(T^{3/2}),\quad\int_0^TE_2^2(t)\d t
\;=\; O(T^2\log^{C_2}T),\leqno(5.3)
$$
with effective constants $C_1,\,C_2 > 0$ (the values $C_1 = 8, C_2 = 22$ are
worked out in [43]). The above results were proved by
Y. Motohashi and the author: (5.2) and the first bound in (5.3) in
[21],[36] and the second upper bound in (5.3) in [35]. The $\Omega$--result
in (5.2)
was improved to $E_2(T) = \Omega_\pm(T^{1/2})$ by Y. Motohashi [42]. Recently
the author [23] made further progress in this problem by proving the following
quantitative omega-result: there exist two constants $A >0,\,B > 1$ such that
for $T \ge T_0 > 0$ every interval $[T,\,BT]$ contains points $T_1,T_2$ for
which
$$
E_2(T_1) \;>\; AT_1^{1/2},\quad E_2(T_2) \;<\; -AT_2^{1/2}.\leqno(5.4)
$$
Very likely this integral is $\sim CT^2$ for some $C > 0$ as $T\to\infty$.
The latest result, proved by the author [32], complements the upper
bound in (5.3) for the mean square, and says  that
$$
\int_0^T E_2^2(t)\d t \;\gg\; T^2.
$$
Very likely this integral is $\sim CT^2$ for some $C > 0$ as $T\to\infty$.
In concluding, let it be mentioned that
the sixth moment was investigated in [29], where it was shown that
$$
\int_0^T|\zt|^6\d t  \ll_\e T^{1+\e}
$$
does hold if a certain conjecture involving the so-called ternary additive
divisor problem is true. On the other hand, it is known that (see Ch. 9 of
[16]) that
$$
\int_0^T|\zt|^{2k}\d t \gg_k T(\log T)^{k^2}\qquad(k \in \NN).
$$

\smallskip
\centerline{\bb References}
\smallskip
{\hh
\item{[1]} F.V. Atkinson, The mean value of the Riemann
zeta-function,  Acta Math., { 81} (1949), 353--376.
\item{[2]} R.C. Baker and G. Harman,  The difference between consecutive
primes, London Math. Soc. (3)72(1996), 261-280.

\item{[3]} J.L. Hafner and A. Ivi\'c, On the mean square of the
Riemann zeta-function on the critical line.  J. Number Theory, 
 32 (1989), 151--191.
\item{[4]}  D.R. Heath-Brown, The fourth moment of the Riemann 
zeta-function,
 Proc. London Math. Soc.  (3)38(1979), 385-422.
\item{[5]} D.R. Heath-Brown,  The distribution and moments of the error
term in the Dirichlet divisor problem, Acta Arith.  60(1992), 389-415.
\item {[6]} M.N. Huxley, On the difference between consecutive primes,
Invent. Math. 15(1972), 155-164.
\item {[7]} M.N. Huxley, A note on exponential sums with a
difference,  Bull. London Math. Soc.,  29(1994), 325--327.
\item{[8]} M.N. Huxley, Area, Lattice Points and Exponential Sums,
Clarendon Press, Oxford, 1996.
\item{[9]}  A.E. Ingham, Mean-value theorems in the theory of the Riemann
zeta-function,  Proc. London Math. Soc.  (2)27(1926), 273-300.
\item{[10]}  A.E. Ingham, On the estimation of $N(\s,T)$, Quart. J.
Math. (Oxford), 11(1940),  291-292.

\item{[11]} A. Ivi\'c, A note on the zero-density estimates for 
the zeta function, Archiv  Math. (Basel-Stuttgart) 33(1979), 155-164.

\item{[12]} A. Ivi\'c,  Large values of the error term in the divisor
problem, Invent. Math.  71(1983), 513-520.

\item{[13]} A. Ivi\'c, Topics in recent zeta-function theory, Publications 
Math\'ematiques d'Orsay, Universit\'e Paris-Sud, Orsay 1983, 272 pp.

\item{[14]} A. Ivi\'c, Exponent pairs and power moments of the zeta-function,
  Proc. 
    Budapest Conf. Number Theory July 1981, Coll. Math. Soc. J. Bolyai 34,
    North-Holland, Amsterdam 1984, 749-768.

\item{[15]} A. Ivi\'c, A zero-density theorem for the Riemann zeta-function, 
Proceedings Moscow 
    Conference in Number Theory, September 1981, Trudy Mat. Inst. AN SSSR
    {163} (1984), 85-89.

\item{[16]} A. Ivi\'c, The Riemann zeta-function, John Wiley and Sons, 
New York
    1985, XVI + 517 pp.

\item{[17]}A. Ivi\'c,  On consecutive zeros of the Riemann
zeta-function on the critical line.  S\'emin. de Th\'eorie des
Nombres, Universit\'e de Bordeaux 1986/87, Expos\'e no.  29, 14 pp.

\item{[18]} A. Ivi\'c, On a problem connected with zeros of
$\zeta(s)$ on the critical line.  Monatshefte Math.,  104 
(1987), 17--27.

\item{[19]} A. Ivi\'c, 
 Les z\'eros de la fonction zeta de Riemann sur la droite critique, Groupe de 
    travail en th\'eorie analytique et \'elementaire des nombres 1986-1987, 
    Universit\'e Paris-Sud, Orsay, 88-01, pp. 47-51.

\item{[20]} A. Ivi\'c,
 Large values of certain number-theoretic error terms, Acta Arith.,
 56(1990), 135-159.

\item{[21]} A. Ivi\'c,  The mean values of the Riemann zeta-function,
LNs  82, Tata Inst. of Fundamental Research, Bombay (distr. by
Springer Verlag, Berlin etc.), 1991.

\item{[22]} A. Ivi\'c,  On the fourth moment of the Riemann
zeta-function, Publs. Inst. Math. (Belgrade)  57(71)(1995), 101-110.

\item{[23]}  A. Ivi\'c, The Mellin transform and the Riemann zeta-function, 
in ``Proceedings of the Conference on Elementary and Analytic Number Theory 
(Vienna, July 18-20, 1996)", Universit\"at Wien \& Universit\"at f\"ur
Bodenkultur, Eds. W.G. Nowak and J. Schoi{\ss}engeier, Vienna 1996, 112-127.

\item{[24]} A. Ivi\'c, On a class of convolution functions connected
with $\zeta(s)$,  Bull.\ CIX Acad.\ Serbe des Sciences
et des Arts, Sci.\ Math.,  20(1995), 29--50.

\item{[25]} A. Ivi\'c,
 On sums of gaps between the zeros of $\zeta(s)$ on the critical line,
      Univ. Beograd. Publ. Elektrotehn. Fak. Ser. Mat. 6(1995), 55-62.

\item{[26]} A. Ivi\'c, 
On the fourth moment of the Riemann zeta-function, Publications Inst.
      Math. (Belgrade) 57(71)(1995), 101-110.

\item{[27]} A. Ivi\'c, On the distribution of zeros of a class
of convolution functions,  Bull. CXI Acad.\ Serbe des Sciences
et des Arts, Sci. Math.,  21(1996), 61--71.

\item{[28]} A. Ivi\'c,
 Uvod u Analiti\v cku Teoriju Brojeva, Izdava\v cka knji\v zarnica
	Zorana Stojanovi\'ca, Sremski Karlovci -- Novi Sad, 1996, 389 str.

\item{[29]} A. Ivi\'c,
 On the ternary additive divisor problem and the sixth moment of the
zeta-function, ``Sieve Methods, Exponential Sums, and their Applications in
Number Theory" (eds. G.R.H. Greaves, G. Harman, M.N. Huxley), Cambridge
University Press, Cambridge, 1996, 205-243.

\item{[30]} A. Ivi\'c, On some results concerning the Riemann Hypothesis, in 
``Analytic Number Theory", LMS LNS 247, Cambridge University Press, 
Cambridge, 1997, pp. 139-167.

\item{[31]} A. Ivi\'c, 
On some problems involving the mean square of $|\zeta({\txt{1\over2}}
 + it)|$,
Bulletin CXVI de l'Acad\'e\-mie Serbe des Sciences et des Arts - 1998, 
Classe des Sciences math\'ematiques et naturelles, 
Sciences math\'ematiques No. 23, pp. 71-76.

\item{[32]} A. Ivi\'c,
 On the error term for the fourth moment of the Riemann zeta-function,
Journal London Math. Society 60(2)(1999), 21-32.

\item {[33]} A. Ivi\'c and M. Jutila, Gaps between consecutive
zeros of the Riemann zeta-function, Monatshefte Math., 
105(1988), 59--73.

\item{[34]} A. Ivi\'c and  K. Matsumoto, On the error term in the mean
	square formula for the Riemann zeta-function in the critical strip, 
	Monatshefte Math. 121(1996), 213-229.

\item{[35]} A. Ivi\'c and Y. Motohashi, The mean square of the
error term for the fourth moment of the zeta-function,  Proc. London Math.
Soc. (3)66(1994), 309-329.

\item {[36]}  A. Ivi\'c and Y. Motohashi,  The fourth moment of the
Riemann zeta-function,  J. Number Theory  51(1995), 16-45.

\item{[37]} A. Ivi\'c and H.J.J. te Riele, On the zeros of the error
term for the mean square of $\;\vert\zeta({\txt{1\over 2}}  + it)\vert\;$.
  Math. Comp., 56 No.  193(1991), 303-328.

\item{[38]} M. Jutila, Zero-density estimates for $L$--functions,
Acta Arith. 32(1977), 52-62.

\item{[39]} N. Levinson, More than one third of the zeros of the
zeta-function are on $\s = \hf$, Advances Math. 13(1974), 383-436.

\item {[40]} J. van de Lune, H.J.J. te Riele and D.T. Winter,
On the zeros of the Riemann zeta-function in the critical strip IV,
 Math. Comp.,  46 (1987), 273--308.

\item{[41]} K. Matsumoto, Recent developments in the mean square theory
of the Riemann zeta and other  zeta-functions, Indian. Nat. Acad. Sci.,
to appear.

\item {[42]}  Y. Motohashi, A relation between the Riemann zeta-function
and the hyperbolic Laplacian,  Ann. Sc. Norm. Sup. Pisa, Cl. Sci. IV
ser.  22(1995), 299-313.

\item {[43]}  Y. Motohashi,  Spectral theory of the Riemann
zeta-function,  Cambridge University Press, Cambridge, 1997.

\item{[44]} A.M. Odlyzko, On the distribution of 
spacings between the zeros of the zeta-function,  Math. Comp.,
48 (1987), 273--308.
\item{[45]} A.M. Odlyzko,  Analytic computations in number theory,
 Proc. Symp. Applied Math.,  48(1994), 451--463.
\item{[46]} A.M. Odlyzko,  The  $10^{20}$-th zero of the Riemann
zeta-function and 175 million of its neighbors, to appear.

\item{[47]} E. Preissmann, Sur la moyenne de la fonction z\^eta,
in ``Analytic Number Theory and Related Topics" (ed. K. Nagasaka),
World Scientific, Singapore, 1993, pp. 19-125.

\item {[48]} H.J.J. te Riele and J. van de Lune,
Computational  number theory  at CWI in 1970--1994, CWI Quarterly, 
7(4) (1994), 285--335.
\item{[49]} B. Riemann, \"Uber die Anzahl der Primzahlen unter
einer gegebenen \ Gr\"osse,   Monats.\ Preuss.\ Akad.\ Wiss., 
(1859--1860), 671--680.
\item {[50]} A. Selberg, On the zeros of  Riemann's zeta-function,
 Skr.\ Norske Vid.\ Akad.\ Oslo,  10(1942), 1--59.
\item {[51]} E.C. Titchmarsh,  The theory of the Riemann zeta-function,
 Clarendon Press, Oxford, 1951.
\item{[52]} K.--M. Tsang,  Higher power moments of $\D(x), E(t)$ and
$P(x)$, Proc. London Math. Soc. (3)65(1992), 65-84.

}

\bigskip\bigskip
\no{\sevenrm
\no Aleksandar Ivi\'c

\no Katedra Matematike RGF-a Universiteta u Beogradu

\no \DJ u\v sina 7, 11000 Beograd, Serbia (Yugoslavia)}

\no \sevenbf e-mail: aivic@matf.bg.ac.yu, eivica@ubbg.etf.bg.ac.yu

\bye